\documentclass{amsart}
\usepackage[procnames]{listings}
\usepackage{pdfpages}
\usepackage[procnames]{listings}
\usepackage{color}
\usepackage{hyperref}
\numberwithin{equation}{section}
\usepackage{cite}
\usepackage{relsize}
\usepackage{amssymb}
\numberwithin{equation}{section}
\usepackage{amsmath}
\usepackage{mathtools}
\usepackage{amssymb}
\usepackage{tikz}
\usepackage{tikz-cd}
\usetikzlibrary{matrix,arrows,decorations.pathmorphing}
\title{ Bounds on the  multiplicity of the Hecke eigenvalues}
\author{Naser T. Sardari}
\email{ntalebiz@math.wisc.edu}
\address{Department of Mathematics, UW-Madison, Madison, WI 53706}
\date{\today}
\usepackage{amsmath}

\usetikzlibrary{positioning}

	\newtheorem{thm}{Theorem}[section]
	
	\newtheorem{prop}[thm]{Proposition}
	
        \newtheorem{conj}[thm]{Conjecture}

	\newtheorem{lem}[thm]{Lemma}

	\newtheorem{cor}[thm]{Corollary}
	
	\theoremstyle{defi}

	\theoremstyle{pf}

	\numberwithin{equation}{section}

\begin {document}
\maketitle
\begin{abstract}Fix an integer $N$ and a prime $p\nmid N$ where $p\geq 5$. 
We show that the number of newforms $f$ (up to a scalar multiple) of level $N$ and even weight $k$ such that $\mathcal{T}_p(f)=0$  is bounded independently of $k$, where $\mathcal{T}_p(f)$ is the Hecke operator.
 \end{abstract}
\section{Introduction}
\subsection{Motivation}
We begin by describing the  \textit{spectral multiplicity hypothesis} introduced in the work of Phillips and Sarnak \cite{Phillips},  and motivating our results from the analytic point of view. Let $\mathbb{H}$ be the upper half-plane and consider the Hilbert space $L^{2}(\mathrm{SL}_2(\mathbb{Z})\backslash \mathbb{H})$ equipped with the Laplace-Belteram operator $\Delta$ induced from the hyperbolic metric $\frac{1}{y^2} (dx^2+dy^2)$ on $\mathbb{H}.$ Selberg established the spectral  decomposition  of $\Delta$ acting on $L^{2}(\mathrm{SL}_2(\mathbb{Z})\backslash \mathbb{H})$:
$$
 L^{2}(\mathrm{SL}_2(\mathbb{Z})\backslash \mathbb{H})= L^{2}_{\text{cusp}}(\mathrm{SL}_2(\mathbb{Z})\backslash \mathbb{H})+\text{Eis}
$$
where $\text{Eis}$ is the contribution of the Eisenstein series $\big($form the continous part of the spectrum filling out $[1/4,\infty)$$\big)$ and its possible residues  (the constant function which is the residue of the pole at $s=1$), and $ L^{2}_{\text{cusp}}(\mathrm{SL}_2(\mathbb{Z})\backslash \mathbb{H})$ is the space of Maass cusp forms; see  \cite{Selberg}. It is widely believed that the entire cusp spectrum is simple; and this is very difficult to prove. The spectral multiplicity hypothesis is the assumption that a positive density of these cusp eigenvalues have a uniformly bounded multiplicity.

 This assumption and its  variations  are essential in the deformation theory of hyperbolic surfaces that is initiated by Phillips, Sarnak and later developed by Wolpert.  It is proved that, under this assumption, a generic hyperbolic curve with  cusp, has only finitely many cusp eigenvalues unlike $
 L^{2}(\mathrm{SL}_2(\mathbb{Z})\backslash \mathbb{H})$ or the arithmetic lattices  where the cuspidal  spectrum dominates the continuous spectrum.

The best bounds on the multiplicity of the eigenvalues are far from this conjecture. More precisely,  let $m(t)$ be the multiplicity of the eigenvalue $1/4+t^2$ on $L^{2}(\mathrm{SL}_2(\mathbb{Z})\backslash \mathbb{H})$.  The best known upper  bound is 
\begin{equation}\label{mult}
m(t) \leq \frac{\pi t}{24\log(t)}(1+o(1)),
\end{equation}
where  $o(1)\to 0$ as $t\to \infty$; see  \cite{Sarnak3}.
 As pointed out by Sarnak  \cite [Page 2]{Sarnak3}, even showing $m(t)=o(t/\log(t))$ is very difficult. The problem is that the Laplacian  eigenvalues for large $t$ are not isolated (the expected consecutive distance is $1/t$) and this makes the problem inaccessible by analytic methods.  

Bounding the multiplicity of the laplacian eigenvalues also appears in the quantum uniques ergodicity.  It is known that the existence of a sequence of  eigenvalues with large multiplicity violates the QUE conjecture \cite{Sarnak2}. It is expected that for every compact hyperbolic manifold $m(t)=t^{\epsilon},$ for any $\epsilon>0$ and again the best known bound is $m(t)=O(t^{d-1}/\log(t))$ \cite{Pierre}, where $d$ is the dimension of the manifold.  

In this paper, we consider the problem of bounding the multiplicity of Hecke operator's eigenvalues for the family of  the holomorphic modular forms with the  fixed level and varying weight.

Let $
\Gamma_0(N):=\left\{\begin{bmatrix}a & b \\ c &d  \end{bmatrix}\in \mathrm{SL}_2(\mathbb{Z}): a, b, c,d \in \mathbb{Z}, N|c   \right\}
$  be the Hecke congruence subgroup of level $N.$ Let  $S_{k,\chi}(N)$ be the space of holomorphic cusp forms of even weight $k\in \mathbb{Z}$, level $N$, and  nebentypus character $\chi$. It is the space of the holomorphic functions $f$ on the upper half plane $\mathbb{H}$ such that
\begin{equation}\label{transform}
f\left(\frac{az+b}{cz+d}\right)=(cz+d)^k \chi(d) f(z)
\end{equation}
for every $\begin{bmatrix} a &b \\ c & d \end{bmatrix} \in \Gamma_0(N)$, and $f$ converges to zero as it approaches each cusp (we have finitely many cusps for $\Gamma_0(N)$ that are associated to the orbits of $\Gamma_0(N)$ acting by M\"obius transformations on $\mathbb{P}^1(\mathbb{Q})$).   It is well-known that $S_{k,\chi}(N)$ is a finite dimensional vector space over $\mathbb{C}$, and is equipped with the Petersson inner product $\langle f, g \rangle:=\int_{\Gamma_0(N) \backslash \mathbb{H}} f(z)\bar{g}(z) y^k dx dy/y^2$ which makes it into a Hilbert space. Assume that $p$ is a fixed  prime number where $p\nmid N$. Then  one can define a self-adjoint  Hecke operator $\mathcal{T}_p$ on $S_k\big(N\big)$:
\begin{equation}\label{Hecke}
\mathcal{T}_p(f)(z):=\sum_{n=1}^{\infty} a_{np}e(nz)+p^{k-1}\chi(p)\sum_{n=1}^{\infty} a_{n}e(pnz),
\end{equation}
where $f(z)=\sum_{n=1}^{\infty} a_{n}e(nz)$ is the Fourier expansion of $f$ at the cusp $\infty$. In particular, if $f$ is an eigenfunction of $\mathcal{T}_p$ with eigenvalue $\lambda_p(f)$ then
$a_p=a_1\lambda_p(f).$ By Deligne's result \cite{Deligne,Pierre} the Ramanujan-Petersson conjecture  holds for $f$ and we have $|\lambda_p(f)|\leq 2p^{\frac{k-1}{2}}$.
Under Langlands' philosophy,  the Hecke operator $\mathcal{T}_p$ is the $p$-adic analogue of the Laplace operator \big(the eigenvalues of $\mathcal{T}_p$ determine the Satake parameters of the associated local representation $\pi_{p}$ of $GL_2(\mathbb{Q}_p)$ just as the Laplace eigenvalue of the Maass form determines  the associated local representation $\pi_{\infty}$ of  $GL_2(\mathbb{R})$\big).

Let $m_p(\lambda,k,\chi,N)$ be the multiplicity of $\lambda$ as an eigenvalue of $\mathcal{T}_p$ acting on $ S_{k,\chi}(N)$.  We fix the level $N$ in our paper and  the implicit constants in the $O$ notations depend on  $N.$  The trivial upper bound is $m_p(\lambda,k,\chi,N) =O(k)$ which is the dimension of $S_{k,\chi}(N).$ By Eichler-Selberg trace formula  Serre proved that $m_p(\lambda,k,\chi,N) =o(k)$ \cite{Serre0} and it is not hard to derive from his work that  $m_p(\lambda,k,\chi,N) =O(k/\log(k))$; see \cite{Murty}. This is the analogue of \eqref{mult} in the weight aspect. Serre derived a number of striking consequences from this  bound on the multiplicities. For example, he proved that if $(N_{i})$ is any sequence of positive integers and if $(d_i)$ denotes the maximum of the dimensions of the simple abelian variety quotients of $Jac(X_0(N_{i}))$, then $d_i\to \infty$ as $N_i\to \infty$. In particular, there are only finitely many positive integers $N$ for which $Jac(X_0(N_{i}))$ is isogenous to a product of elliptic curves.

For $\lambda\neq 0$, Frank Calegari in his blog post \cite{Calegari} proved  that
\begin{equation}\label{nonzero}
m_p(\lambda,k,\chi,N)=O(v_p(\lambda)^2),
\end{equation}
where the implicit constant in $O(.)$ is independent of $k$ and only depends on fixed numbers $N$ and  $p.$ We briefly explain his method. Let $f$ be a modular form with $p$-th Hecke eigenvalue $\lambda\neq 0.$ Then the  slope of the modular form $f$ is finite. Since the level is fixed, there are only finitely many Coleman families of modular forms \cite{col1, col2} that cover these modular forms. \eqref{nonzero} follows from Wan's explicit quadratic bound \cite{Wan} on Gouvea and Mazur conjecture \cite{Mazur}. This method does not work for bounding $m_p\left(0,k,\chi,N\right).$
\\

 In this paper we use methods in the deformation theory of Galois representations and the Taylor-Wiles method  to give a quantitative bound on  $m_p\left(0,k,\chi,N\right)$ which are independent of $k$. 
\begin{thm}\label{mainthm} Suppose that $p\geq 5$ is prime and $k>0$ is an even integer.  
We have 
\begin{equation}
m_p\left(0,k,\chi,N\right)\leq m_p(0,k^{\prime},\chi,N),
\end{equation}
for some $k^{\prime}\leq \frac{p+3}{2}$ which is explicitly given in Theorem~\ref{serreconj} in terms of $k.$
\end{thm}
 Finally, we construct a family of newforms $f \in S_{k,\chi}(N)$ with $a_p(f)=0$ for any $k\geq 1$ and prime $p.$
Let $D<0$ be a square-free integer, $\mathcal{O}_{\sqrt{D}}$ be   the ring of integers of $\mathbb{Q}(\sqrt{D}),$  and $\chi_D(n):= \left(\dfrac{D}{n}  \right)$ be the quadratic character associated to $D$. Suppose that $\xi$ is an algebraic Hecke character of $\mathcal{O}_{\sqrt{D}}$ with conductor $\mathfrak{m}\subset \mathcal{O}_{\sqrt{D}}$ and weight $k,$ where $\mathfrak{m}$ is prime to $D.$ Let $N:=|D|\text{Nr}(\mathfrak{m}),$ where $\text{Nr}(\mathfrak{m}):= \# \left(\mathcal{O}_{\sqrt{D}}/ \mathfrak{m} \right).$
Let
\[
f_{\xi}(z):=\sum_{\mathfrak{a}\subset \mathcal{O}_{\sqrt{D}}} \text{Nr}(\mathfrak{a})^{k/2} \xi(\mathfrak{a}) e(z\text{Nr}(\mathfrak{a})),
\]
where $\sum_{\mathfrak{a}\subset \mathcal{O}_{\sqrt{D}}} $ is over the ideals $\mathfrak{a}$ of $ \mathcal{O}_{\sqrt{D}}.$ By Theorem~\cite[Theorem 12.5]{Iwaniec}, we have $f_k \in S_{k+1,\chi}(N)$ is a newform, where $\chi(n):= \chi_D(n)\xi((n))$ for $n\in \mathbb{Z}.$ We note that if $p$ does not split in $\mathbb{Q}(\sqrt{D}),$ then $a_p(f_k)=0.$ We say $f_{\xi}$ is a CM newform. Let $m_{p,\text{CM}}(0,k,\chi,N)$ denote the number of CM  newforms $f \in S_{k,\chi}(N)$ with $a_p(f)=0.$ Finally, we make the following conjecture.
\begin{conj}
Let $p$ be a prime number and $N \geq 0$ be an integer, where $\gcd(p,N)=1$. We have 
\[
m_p\left(0,k,\chi,N\right)= m_{p,\text{CM}}(0,k,\chi,N),
\] 
for large enough weight $k$ which only  depends on $p$ and $N.$
\end{conj}
 We give an outline of our proof for bounding $m_p\left(0,k,\chi,N\right)$.  
\subsection{Outline of the proof}
 We fix a modular absolutely irreducible  residual  representation   $\bar{\rho} : G_{\mathbb{Q}} \to GL_2(\mathbb{F}_{p^n})$ with prime to $p$ conductor dividing $N$ and consider the problem only for $f\in S_{k,\chi}(N)$ whose associated mod $p$ representation (which is well-defined since $\bar{\rho}$ is absolutely irreducible) is isomorphic to $\bar{\rho}$. Note that by  (part of) Serre's conjecture \cite{strongserre} there are only a bounded number (only depends on prime $p$ and $N$ and not on $k$) of such  $ \bar{\rho}$; see Theorem~\ref{serreconj}. Moreover,  since  $a_p(f)=0$,  the restriction of the $p$-adic representation $r_f:G_{\mathbb{Q}}\to GL_2(\bar{\mathbb{Q}}_p) $ associated to $f$   to the decomposition group at $p$ is a dihedral representation (induced from a Lubin-Tate character of $G_{\mathbb{Q}_{p^2}}$, where  $\mathbb{Q}_{p^2}$ is the unramified degree 2 extension of $\mathbb{Q}_p$ with ring of integers $\mathbb{Z}_{p^2}$); see Theorem~\ref{dihedthm}. Moreover,  it has irreducible reduction if $k$ is not congruent to $1$ mod $p+1$ and in particular has irreducible reduction when $p>2$ and $k$ is even; see Theorem~\ref{dihedthm}. So, we study the deformation ring of deformations of $\bar\rho$ unramified outside $Np$ and whose restriction to the inertia subgroup of the decomposition group at $p$ is dihedral and at other places has a fixed inertial type. We construct a deformation ring $R_{D}$ such that any $r_f$ factors through it.    $R_{D}$ will be a $\Lambda_p$-algebra (the deformation ring associated to the characters of  $G_{\mathbb{Q}_{p^2}}$). By the class field theory  $\Lambda_p$ is isomorphic to the power series ring $\mathbb{Z}_{p^n}[[S_1,S_2,S_3]]$ (at least if $p>2$). To an integer $k>1$ there is a dimension 1 prime ideal $P_k$ of $\Lambda_p$ corresponding to the character of $1+p\mathbb{Z}_{p^2}$ sending $a$ to $a^{k-1}$ and the unramified  character sending  Frobenius element to $\sqrt{-1};$ see Theorem~\ref{dihedthm}.  If $f$ is an newform of  weight $k$ then $r_f$ corresponds to a point of $(R_{D}/P_kR_{D})^{\text{red}}[1/p]$. So it would suffice to show that $\dim_{\mathbb{Q}_{p^n}} (R_{D}/P_kR_{D})^{\text{red}}[1/p]$ is bounded independent of $k$. Equally it would suffice to show that $R_{D}$ is finitely generated as a $\Lambda_p$-module. By the topological form of Nakayama's lemma (see \cite[Exercise 7.2]{Eisenbud}) it suffices to show that for one prime ideal $P$, $R_{D}/PR_{D}$ is finitely generated as a $\mathbb{Z}_{p^n}$-module. If $p\geq 5$ then we show that there exists  a prime ideal $P_{min}$ of $\Lambda_p$ so that the Hodge-Tate weights are of moderate Hodge-Tate type \cite{FM} (in the Fontaine-Laffaille range, we use strong Serre's conjecture for this step). By a modularity result in the work of Khare and Wintenberger, we show that there exists a finite $\mathbb{Z}_{p^n}$ algebra $ \bar{R}_S^{\psi}$ with a surjection to $R_{D}/P_{min}R_{D}.$ We show that $ \bar{R}_S^{\psi}$   is isomorphic to  an appropriate  Hecke algebra which gives explicit bound on the number of generators of  $ \bar{R}_S^{\psi}$ and $R_{D}/P_{min}R_{D}$ as  finite  $\mathbb{Z}_{p^n}$-modules. 
 \subsection{Notations} We assume that $N$ is a fixed integer and  $p$ is a fixed prime number where $\gcd(p,N)=1.$
We denote the space of modular forms of level $N$ and weight $k$ with  nebentypus character $\chi:\big(\mathbb{Z}/N\mathbb{Z}\big)\to \mathbb{C}^*$ by  $S_{k,\chi}(N)$ which is a finite dimensional vector space of dimension $\frac{\alpha(N)k}{12}+O(N^{1/2+\epsilon}),$ where $\alpha(N)=N\prod_{p|N}(1+1/p).$ Define  $S_k(N):= \oplus_{\chi} S_{k,\chi}(N)$. For a newform $f\in S_k(N)$, we write  $f(z)=\sum_{n=1}^{\infty} a_{n}e(nz)$ for the Fourier expansion of $f$ at cusp $\infty$  and wirte $a_p(f)$ for the $p$-th Fourier coefficient of $f.$  Let  $G_{\mathbb{Q}}:=\text{Gal} (\bar{\mathbb{Q}}/{\mathbb{Q}})$,  $G_{\mathbb{Q}_p}:=\text{Gal} (\bar{\mathbb{Q}}_p/{\mathbb{Q}_p})$, $I_p$ be  the inertia subgroup at $p$  and $\mathbb{Q}_{p^n}$ be the unique unramified extension of degree $n$ of $\mathbb{Q}_p$ with ring of integer $\mathbb{Z}_{p^n}$ and residue field $\mathbb{F}_{p^n}.$ We write $\varepsilon:G_{\mathbb{Q}} \to \mathbb{Z}_p^*$ for the cyclotomic character. We denote the local reciprocity map of the local field $\mathbb{Q}_{p^n}$ by $rec_{p^n}: \mathbb{Q}_{p^n}^{*} \to G_{\mathbb{Q}_{p^n}}^{ab}$ which is well-defined up to an embedding $\mathbb{Q}_{p^n}\hookrightarrow \bar{\mathbb{Q}}_{p}$.  We fix an embedding $i: \bar{\mathbb{Q}}\to \bar{\mathbb{Q}}_p$ which  defines a $p$-adic valuation $v_p.$ Given a  newform $f\in S_{k,\chi}(N) $, we denote the associated  $p$-adic representation by $\rho_f: G_{\mathbb{Q}}\to GL_2(\bar{\mathbb{Q}}_p)$,  its restriction  to the decomposition group at $p$ by  $ \rho_{f,p}: G_{\mathbb{Q}_p} \to GL_2(\bar{\mathbb{Q}}_p) $  and its residual representation if it is irreducible  by $\bar{\rho}_f: G_{\mathbb{Q}}\to GL_2(\bar{\mathbb{F}}_p)$ which are well-defined up to conjugation; see \cite{Deligne}. Define the two Lubin-Tate characters $\varepsilon_2, \varepsilon^{\prime}_2 : G_{\mathbb{Q}_{p^2}} \to \mathbb{Z}_{p^2}^*$ by requiring  $\varepsilon_2,\varepsilon_2^{\prime} \circ rec_{p^2} (p)=1$ and $\varepsilon_2,\varepsilon_2^{\prime} \circ rec_{p^2}|_{\mathbb{Z}_{p^2}^*} = \mathbb{Z}_{p^2}\hookrightarrow \bar{\mathbb{Q}}_p$ given by two natural embedding of $ \mathbb{Z}_{p^2}\hookrightarrow \bar{\mathbb{Q}}_p$.    Let $\mu_{a}:G_{\mathbb{Q}_p}\to \bar{\mathbb{Q}}_l$ be the unique unramified character such that  $\mu_a(\text{Frob}_p)=a.$

\section{Modular forms with $a_p(f)=0$ }
Let $f\in S_{k,\chi}(N)$ be a newform with $a_p(f)=0$,  where $\gcd(p,N)=1.$  We cite a theorem which shows  that the restriction of the Galois representation $\rho_f$ to the decomposition group at $p$ is a dihedral representation. 

\begin{thm}[Breuil] \label{dihedthm}
We have 
\begin{equation}
\rho_{f|G_{\mathbb{Q}_p}}= \text{Ind}_{G_{\mathbb{Q}_p^2}}^{G_{\mathbb{Q}_p}} \varepsilon_2^{k-1} \otimes \mu_{\sqrt{-1}}: G_{\mathbb{Q}_p} \to GL_2(\mathbb{Q}_p).
\end{equation}
Moreover,  $ \text{Ind}_{G_{\mathbb{Q}_p^2}}^{G_{\mathbb{Q}_p}} \varepsilon_2^{k-1} \otimes \mu_{\sqrt{-1}}$ is a crystalline representation of $G_{\mathbb{Q}_p}$ for any integer $k\geq 2$, and its mod $p$ reduction is absolutely irreducible if $k$ is not congruent to $1$ mod $p+1$. In particular $\bar{\rho}_f$ is absolutely  irreducible when $p>2$ and $k$ is even.
\end{thm}
\begin{proof}
See~\cite[Proposition 3.1.2]{Breuil}.
\end{proof}
\subsection{Mod $p$ representations}
We cite a theorem which shows that every modular form with $p$-th Fourier coefficient zero  is congruent up to a twist with a modular form of weight at most $ (p+3)/2$ mod $p.$ Note that for $p\geq 5,$ $(p+3)/2\leq p-1$ which is in the Fontaine-Laffaille range. 

\begin{thm}[Ash and Stevens - Edixhoven]\label{serreconj}
Let $h\in S_{k,\chi}(N)$  be a newform then there exist integers 
$i$ and $k^{\prime}$ with $0 \leq i\leq p-1, k^{\prime} \leq p+1$ and a modular form  $g\in S_{k^{\prime},\chi}(N)$ such that $\bar{\rho}_{h} \sim \varepsilon^i \otimes \bar{\rho}_g.$ Moreover, assume that $a_p(h)=0$, $p$ is an odd prime number and $k$ is an even integer. Let $0\leq a< b \leq p-1$ be uniquely determined by the congruence condition  $k-1=b+pa \text{ or } a+bp \mod p^2-1$. Then $\bar{\rho}_{h}$ is absolutely irreducible and   the only possible values for $(i,k^{\prime})$ which all occur are:
$$
\begin{cases}
(a,1+b-a), (b-1,p+2+a-b)  &\text{ if } b-a\neq 1
\\
(a,2)  &\text{ if } b-a=1.
\end{cases}
$$
Note that we can choose $k^{\prime}\leq (p+3)/2$ and since $p\geq 5$ then $k^{\prime}$ is in the Fontaine-Laffaille range. 
\end{thm}
\begin{proof}
 See \cite[Theorem 3.4]{Edi}  and  \cite[Theorem 4.5]{Edi}.
\end{proof}

\section{Deformation rings of Galois representations}\label{cryst}
In this section,  we define several deformation rings which are complete local noetherian rings with fixed finite  residue field of characteristic $p$. We cite a result of  Khare and Wintenberger~\cite[Theorem 10.1]{Khare} on finiteness of deformation rings with some conditions. We keep our notations consistent with their work in this section and work with framed deformation rings in order to cite their results. We refer the reader to the notes of Gebhard \cite{Gebhard} and the work of Mazur \cite{Mazur} for an introduction on the deformations of Galois representations.

\subsection{Galois deformation rings with local conditions}
Let $S$ be the union of the set of primes in $pN.$ Let $\bar{\rho}:G_{\mathbb{Q}}\to GL_2(\mathbb{F}_{p^{n}})$ be a mod $p$ representation of a newform $h\in S_{k,\chi}(N)$ with $a_p(h)=0$ for some weight $k$. By Theorem~\ref{serreconj}, $\bar{\rho}$ is an absolutely irreducible   representation  of conductor dividing $N$ and the number of them are bounded. We fix $\bar{\rho}$ in this section and suppose that $p\geq 5$. Let $\mathcal{C}_{\mathbb{Z}_{p^n}}$ be  the category of the complete local noetherian $\mathbb{Z}_{p^n}$-algebras with  residue field $\mathbb{F}_{p^n}$.

  For $v\in S$, let $R_v^{\square}$ denote the local universal framed deformation ring of $\bar{\rho}|_{G_{\mathbb{Q}_v}}$ in $\mathcal{C}_{\mathbb{Z}_{p^n}}.$
Let $\psi: G_{\mathbb{Q}} \to \bar{\mathbb{Q}}_{p}^{*}$  be a continuous homomorphism such that $\psi= det (\bar{\rho}) \mod p.$ We write $\psi_v: G_{\mathbb{Q}_v} \to \bar{\mathbb{Q}}_{p}^{*}$ for the restriction of $\psi$ to the local Galois groups.  Let $R_v^{\square,\psi}$ be  the local framed deformation ring of $\bar{\rho}|_{G_{\mathbb{Q}_v}}$ with the determinant condition $det=\psi.$   For $v\in S$, we define a certain local deformation rings which are quotient of  $R_v^{\square}.$ We denote them by  $\bar{R}_p^{\square,\psi}$ if $v=p,$ and $R_{v,\tau_v}^{\square,\psi}$ for an inertial type $\tau_v.$   These classify  a set of  prescribed conditions $X_v$, including fixed determinant $\det(\rho)=\psi $.  We also verify the required properties  of $\bar{R}_p^{\square,\psi}$ and $R_{v,\tau_v}^{\square,\psi}$ stated  in \cite[Theorem 3.1]{Khare} which are necessary  to apply  the finiteness theorem in~\cite[Theorem 10.1]{Khare}. 
\subsubsection{ $\bar{R}_p^{\square,\psi_p}$: low weight crystalline case}\label{lowweight}
We cite a theorem which is essentially due to Ramakrishna \cite[Theorem 4.2]{Ramakrishna}. We need a version of this theorem  for the framed deformation rings with more general  determinant conditions.  We cite a version of this theorem from the work of Fontaine and Mazur  \cite[Theorem B2]{FM}.  
\begin{thm}\label{rama} Recall $p,$ $\psi_p,$ $\bar{\rho}$  and the distinct pair of integers  $(i,k^{\prime})$ associated to $\bar{\rho}$ in Theroem~\ref{serreconj} with $k^{\prime}\leq \frac{p+3}{2}.$ 
There is a unique reduced, $p$-torsion free quotient $\bar{R}_p^{\square,\psi_p}$
of $R_p^{\square}$ with the property that a continuous homomorphism $\alpha: R_v^{\square} \to A $ for $A\in\mathcal{C}_{\mathbb{Z}_{p^n}},$ factors through $\bar{R}_p^{\square,\psi_p}$ if and only if
$\alpha \circ \rho^{\text{univ}}$ is crystalline,    $HT(\alpha \circ \rho^{\text{univ}})=(i,i+k^{\prime}-1)$ and $det(\alpha \circ \rho^{\text{univ}})=\psi_{p}$, where $HT(\alpha \circ \rho^{\text{univ}})$ is the Hodge-Tate weights of $\alpha \circ \rho^{\text{univ}}.$ Furthermore, $\bar{R}_p^{\square,\psi_p}$ is isomorphic to a power series ring in $4$ variables over $\mathbb{Z}_p$.
\end{thm}
\begin{proof}
The proof follows from~ \cite[Theorem B2]{FM} and the fact that $\frac{p+3}{2}\leq p-1$ for $p\geq 5.$
\end{proof}
\subsubsection{$R_{v,\tau_v}^{\square,\psi_v}$ for $v\neq p$ : Deformations of fixed inertial type }\label{type}
We follow closely \cite[Section 3.30]{Gee}.
Given a local representation $\rho: G_{\mathbb{Q}_v} \to GL_2(A),$ where $A\in \mathcal{C}_{\mathbb{Z}_{p^n}},$ there is a Weil-Deligne representation $WD(\rho)$ associated to $\rho.$    If $WD=(r,N^{\prime})$ is a Weil-Deligne representation, then we call $(r|_{I_{\mathbb{Q}_v}},N^{\prime})$ an inertial WD-type, where $r|_{I_{\mathbb{Q}_v}}$ is the restriction of $r$ to the inertia subgroup $I_{\mathbb{Q}_v}$. In particular, if $\rho$ is unramified then its inertial WD-type is $(Id,0).$   We cite the following theorem from \cite[Theorem 3.31]{Gee} and \cite[Section 5.E]{Shotton}.
\begin{thm}\label{Shotton}  Let $\tau_v$ be an inertial  WD-type of conductor $|v|^{\text{ord}_v(N)}$.
There is a unique reduced, p-torsion free quotient $R_{v,\tau_v}^{\square,\psi_v}$ of $R_v^{\square}$ with the property that a continuous homomorphism $\alpha:R_v^{\square}\to A$ for $A\in \mathcal{C}_{\mathbb{Z}_{p^n}}$ factors through $R_{v,\tau_v}^{\square,\psi_v}$ if and only if $\alpha \circ\rho^{\text{univ}}$ has inertial WD-type $\tau$ and $det(\alpha \circ \rho^{\text{univ}})=\psi_{p}.$ Moreover, for all but finitely many $\tau$, we have $R_{v,\tau_v}^{\square,\psi_v}=0$  and if $R_{v,\tau_v}^{\square,\psi_v}$ is nonzero then it is
Cohen-Macaulay, has relative dimension 3 over $\mathbb{Z}_{p^n}$, and 
the generic fibre $R_{v,\tau_v}^{\square,\psi_v}[1/p]$ is irreducible and formally smooth.
\end{thm} 
\subsection{ Finiteness of deformation rings}
We follow the same notations as in \cite[Section 10]{Khare}. Let $R_S^{\square,\psi}$ be the global framed deformation ring of deformations of $\bar{\rho}:G_{\mathbb{Q}}\to GL_2(\mathbb{F}_{p^n})$  with framings at each $v\in S$ and fixed determinant $\psi$ which are unramified  outside $S$. Recall that $R_v^{\square,\psi_v}$ is the local framed deformation ring of deformations of $\bar{\rho}|_{G_{\mathbb{Q}_v}}$ with fixed determinant $\psi_v$.  Define $R^{\square,loc,\psi}:=\hat{\otimes}_{\mathbb{Z}_{p^n}, v\in S}R_v^{\square,\psi_v}$ which is the completed tensor product of $R_v^{\square,\psi_v}$ over $\mathbb{Z}_{p^n}.$ Similarly, define $\bar{R}^{\square,loc,\psi}:=\bar{R}_p^{\square,\psi_p}\hat{\otimes}_{\mathbb{Z}_{p^n},v\in S\backslash \{p\}}R_{v,\tau_v}^{\square,\psi_v}$, where $R_{v,\tau_v}^{\square,\psi_v}$ were defined in sections \ref{cryst} and \ref{type}. Let 
\begin{equation}\label{defring}\bar{R}_S^{\square,\psi}:=R_S^{\square,\psi} \hat{\otimes}_{R^{\square,loc,\psi}}\bar{R}^{\square,loc,\psi}.\end{equation} 
Let $R_{\{\tau_v,\chi_v\}_{v\in S\backslash \{p\}}}$ be the global deformation ring of deformations of $\bar{\rho}$ which are unramified outside $S$  and have the inertial type and the determinant conditions $\{\tau_v,\chi_v\}$ for $v\in S\backslash \{ p\}.$ Let  and $\bar{R}_S^{\psi}$ the image of $R_{\{\tau_v,\chi_v\}_{v\in S\backslash \{p\}}}$ in $\bar{R}_S^{\square,\psi}$. We cite the following theorem from \cite[Theorem 10.1]{Khare}
\begin{thm}\label{finite}
The ring $\bar{R}_S^{\psi}$ is finite as a $\mathbb{Z}_{p^n}$-module.
\end{thm}

\subsection{ $\bar{R}_S^{\psi}$ is Cohen-Macaulay}
We cite a result  of Snowden~\cite[Proposition 5.0.6]{Snowden} which implies  that  $\bar{R}_S^{\psi}$ is flat over $\mathbb{Z}_{p^n}$ and is Cohen-Macaulay. This result is motivated by  the suggestion in the paragraph before Corollary 4.7 in~\cite{Khare},  This shows that  $\bar{R}_S^{\psi}$ ( rather than $\bar{R}_S^{\psi}[1/p]^{\text{red}}$) is  isomorphic to an appropriate completed Hecke algebra $\mathbb{T}_{m}$. Then  we give explicit bounds on the number of generators of $\bar{R}_S^{\psi}$ as a finite $\mathbb{Z}_{p^n}$ module.   
\begin{prop}[Snowden]\label{snow}
Assume $\bar{R}_S^{\psi}$ is finite over $\mathbb{Z}_{p^n}$ and each $R_{v,\tau_v}^{\square,\psi_v}$ is Cohen-Macaulay and equidimensional of dimension $5$ if $v=p$ and 4 if $v\neq p$. Then $\bar{R}_S^{\psi}$ is flat over $\mathbb{Z}_{p^n}$ and Cohen-Macaulay. 
\end{prop}
\begin{cor}\label{maincor}
 $\bar{R}_S^{\psi}$ is finite and  Cohen-Macaulay over $\mathbb{Z}_{p^n}.$
\end{cor}
\begin{proof}
By Theorem~\ref{finite},  $\bar{R}_S^{\psi}$ is finite over $\mathbb{Z}_p.$ By Theorem~\ref{rama} $\bar{R}_p^{\square,\psi}$ is isomorphic to $\mathbb{Z}_{p^n}[[t_1,\dots,t_4]]$. Hence, it has Krull dimension 5 and it is Cohen-Macaulay (formally smooth over $\mathbb{Z}_{p^n}$). By Theorem~\ref{Shotton},  $R_{v,\tau_v}^{\square,\psi_v}$ is
Cohen-Macaulay, has Krull dimension 4. So, all the conditions of Proposition~\ref{snow} are satisfied and  this concludes our theorem.

\end{proof}

 \section{ Proof of  Theorem~\ref{mainthm}  }
We give a proof of Theorem~\ref{mainthm} in this section.  Let $f\in S_{k,\chi}(N)$ be a  newform with $a_p(f)=0.$ Recall that we fixed an embedding $i:\bar{\mathbb{Q}}\to\bar{\mathbb{Q}}_p.$ By \cite[Theorem 3.1]{Richard}, Theorem~\ref{dihedthm} and Theorem~\ref{serreconj}, there exists a Galois representation $\rho_f:G_{\mathbb{Q}}\to GL_2(\mathcal{O}_f),$ where
$E_f$ is a 
finite extension  of $\mathbb{Q}_p$ with the ring of integers $\mathcal{O}_f$ and residue field $\mathbb{F}_{p^n}$ such that 
\begin{itemize}
\item $\rho_f$ is unramified outside $S$ and  $tr(\rho_f(Frob_v))=i(a_v(f))$ for $v\notin S.$
\item $\rho_{f|G_{\mathbb{Q}_p}}= \text{Ind}_{G_{\mathbb{Q}_p^2}}^{G_{\mathbb{Q}_p}} \varepsilon_2^{k-1} \otimes \mu_{\sqrt{-1}}$, and  $\bar{\rho}_{|G_{\mathbb{Q}_p}}$ is absolutely irreducible, so does $\bar{\rho}_f$.
\item $\bar{\rho}_f$ is one of the finite possible $\bar{\rho}$ listed in Theorem~\ref{serreconj}. 
\item For $v\in S$ and $v\neq p$, $\rho_{f|G_{\mathbb{Q}_v}}$ has an inertial type $\tau_v$ with conductor $|v|^{ord_v(N)},$ and  there are only a finite number of such inertial types $\tau_v.$
\item $det(\rho(f))=i(\chi) \varepsilon^{k-1}.$
\end{itemize}
Let $m_p(0,k,\chi,\bar{\rho},\{\tau_v\}_{v \in S\backslash\{p\}})$ be the number of newforms $f$ with the above properties.  We have 
\[
m_p\left(0,k,\chi,N\right)= \sum_{\bar{\rho},\{\tau_v\}_{v\in S\backslash\{p\}}} m_p\left(0,k,\chi,\bar{\rho},\{\tau_v\}_{v\in S\backslash\{p\}}\right),
\]
where the sum is over  finitely many  possible $\bar{\rho}$ and $\{\tau_v\}_{v\in S\backslash\{p\}}.$  We prove a stronger version of Theorem~\ref{mainthm}.
\begin{thm}\label{strong1.1}
Suppose that $p\geq 5$ be a prime and $k>0$ be an even integer. Let $\bar{\rho},$ $k^{\prime},$ and $\{\tau_v\}_{v\in S\backslash\{p\}}$ as above.  We have 
\[
m_p\left(0,k,\chi,\bar{\rho},\{\tau_v\}_{v\in S\backslash\{p\}}\right)\leq m_p\left(0,k^{\prime},\chi,\bar{\rho},\{\tau_v\}_{v\in S\backslash\{p\}}\right).\]
\end{thm}

For the rest of this section, we fix $\bar{\rho}$ and inertial types $\{\tau_v\}_{v\in S\backslash\{p\}}.$ 
\subsection{The dihedral property at $p$} 
Let $\prod_{p,D}$ (the subscript $D$ is for ``dihedral'') be the maximal profinite quotient group of $G_{\mathbb{Q}_p}$  in which the image of $G_{\mathbb{Q}_p^2}$ is abelian and $\pi_{D,p}:G_{\mathbb{Q}_p} \to\prod_{p,D} $ be the natural projection. Let 
$$
\bar{\rho}_{p,k}=  \text{Ind}_{G_{\mathbb{Q}_p^2}}^{G_{\mathbb{Q}_p}} \bar{\varepsilon}_2^{k-1} \otimes \mu_{\sqrt{-1}}.
$$
By Theorem~\ref{dihedthm}, it follows that $\bar{\rho}_{p,k}$ is absolutely irreducible when $p>2$ and $k$ is even. Since $\bar{\rho}_{p,k}$ is a dihedral representation then there exists  
\[\tilde{\rho}_{p,k}:\prod_{p,D} \to GL_2(\mathbb{F}_{p}) \hookrightarrow GL_2(\mathbb{F}_{p^n})\] such that $\bar{\rho}_{p,k}= \tilde{\rho}_{p,k}\circ \pi_{D,p}.$ Let $W_{\mathbb{Q}_{p^2}}$ be the Weil subgroup of $G_{\mathbb{Q}_{p^2}}$ and $rec:\mathbb{Q}_{p^2}^{\times}\to W_{\mathbb{Q}_{p^2}}$ be the Artin reciprocity isomorphism defined by the local class field theory \cite{Serre}.  Let $\Lambda_{p}$ be  the   universal deformation ring of the deformations of   $\tilde{\rho}_{p,k}$ in $\mathcal{C}_{\mathbb{Z}_{p^n}}$ with residue field $\mathbb{F}_{p^n}$; see \cite[Proposition~1]{Mazur}.
\begin{prop}
$\Lambda_{p}$  is isomorphic to $\mathbb{Z}_{p^n}[[t_1,t_2,t_3]]$. 
\end{prop}
\begin{proof}
By the class field theory, it follows that $\Lambda_{p} \cong \mathbb{Z}_{p^n}[[(1+p\mathbb{Z}_{p^2})\times \hat{\mathbb{Z}}]]$. The theorem follows from the well-known  isomorphism $ \mathbb{Z}_{p^n}[[1+p\mathbb{Z}_{p^2}]] \cong \mathbb{Z}_{p^n}[[t_1,t_2]]$; see \cite[Proposition 11]{Mazur}. 
\end{proof}
We introduce a deformation ring associated to modular forms with $p$-th Fourier coefficient zero, inertial types $\{\tau_v\}_{v\in S\backslash\{p\}}$  and fixed mod $p$ residual representation $\bar{\rho}$. We show that this deformation ring is a finite $\Lambda_p$ module, and the number of its generators gives an upper bound on the multiplicity of  $0$ as an eigenvalue of $\mathcal{T}_p$.
Let $(i,k^{\prime})$ be the pair of integers associated to the weight $k$ and the prime $p$ that are defined in Theorem~\ref{serreconj}. It follows from Theorem~\ref{serreconj}  that
$
\text{Ind}_{G_{\mathbb{Q}_p^2}}^{G_{\mathbb{Q}_p}} \varepsilon_2^{k^{\prime}-1} \otimes \varepsilon^i\mu_{\sqrt{-1}}
$
 is a deformation of $\tilde{\rho}_{p,k}$. Let $r_{k^{\prime},i}: \Lambda_{p} \to \mathbb{Z}_{p^n}$ be the  unique   homomorphism associated to $\text{Ind}_{G_{\mathbb{Q}_p^2}}^{G_{\mathbb{Q}_p}} \varepsilon_2^{k^{\prime}-1} \otimes \varepsilon^i\mu_{\sqrt{-1}}$. Let $P_{min}:=ker(r_{k^{\prime},i})$ which is  a prime ideal of $\Lambda_{p}.$

Let  $R_{\{\tau_v,\chi_v\}_{v\in S\backslash\{p\}}}\in \mathcal{C}_{\mathbb{Z}_{p^n}}$ be the global deformation ring of deformations  $\rho$ of  $\bar{\rho}$ which are unramified outside $S$, and   for $v \in S\backslash\{p\}$  have inertial type $\tau_{v}$ with $det(\rho)_v=\chi_v.$   Let $R_{p}$  be the  universal local defamation ring of  $\bar{\rho}_{p}$; see \cite[Proposition~1]{Mazur}.  By the universal properties  of $R_{\{\tau_v,\chi_v\}_{v\in S\backslash\{p\}}}$, $R_{p},$ and $\Lambda_{p}$, we have the following  natural maps induced from the inclusions $i_p:G_{\mathbb{Q}_p} \to G_{\mathbb{Q}}$ and the projection $\pi_{D,p}:G_{\mathbb{Q}_p} \to\prod_{p,D} $
 \begin{center}
 \begin{tikzcd}[column sep=small]
& R_{p} \arrow{dl}[swap]{\tilde{i}_p} \arrow{dr}{\tilde{\pi}_{D,p}}  & \\
R_{\{\tau_v,\chi_v\}_{v\in S\backslash\{p\}}}  & & \Lambda_{p}
\end{tikzcd}
\end{center}
We define
 \begin{equation}\label{R_D}
R_{D}:=R_{\{\tau_v,\chi_v\}_{v\in S\backslash\{p\}}} \hat{\otimes}_{R_p}  \Lambda_{p}.
\end{equation}
and
$$
\rho_D:= G_{\mathbb{Q}} \to GL_2( R_{\{\tau_v,\chi_v\}_{v\in S\backslash\{p\}}}) \to  GL_2(R_{D}).
$$
Assume that $f\in S_{k,\chi}(N)$ is a newform with $a_p(f)=0,$ inertial types $\{\tau_v \}_{v\in S},$ and  $\bar{\rho}_f=\bar{\rho}$. By the universal property of the tensor product, there is a unique map  $\rho_{f,D}:R_{D}\to \mathcal{O}_f$ such that the following diagram commutes 
\begin{equation}\label{diag}
 \begin{tikzcd}[column sep=small]
& R_{p} \arrow{dl}[swap]{\tilde{i}_p} \arrow{dr}{\tilde{\pi}_{D,p}}  & \\
R_{\{\tau_v,\chi_v\}_{v\in S\backslash\{p\}}}\arrow{dr} \arrow{ddr}[swap]{\rho_f} & & \Lambda_{p}  \arrow{dl}  \arrow{ddl}{\rho_{f,D,p}}
\\
&R_{D} \arrow{d}{\rho_{f,D}}
\\
&\mathcal{O}_f
\end{tikzcd}
\end{equation}
Finally we give a proof of Theorem~\ref{strong1.1}.
\begin{proof}[Proof of Theorem~\ref{mainthm}]  
 Note that $R_{D}$ is a $\Lambda_{p}$ algebra. In what follows, we show that in order to bound $m_p\left(0,k,\chi,\bar{\rho},\{\tau_v\}_{v\in S\backslash\{p\}}\right)$,  it is enough to show that $R_{D}$ is a finite $\Lambda_{p}$ module and it is generated by at most $m_p(0,k^{\prime},\chi,\bar{\rho},\{\tau_v\}_{v\in S\backslash\{p\}})$ elements over $\Lambda_{p}$. Assume that $f\in S_{k,\chi}(N)$ is a newform with $a_p(f)=0,$ inertial types $\{\tau_v \}_{v\in S},$ and  $\bar{\rho}_f=\bar{\rho}$. By the diagram~\eqref{diag}, there is a unique map $\rho_{f,D}: R_{D} \to \mathcal{O}_{f}.$ 
 
Let  $P_k$ be the prime ideal in $\Lambda_p$ associated to the representation $\text{Ind}_{G_{\mathbb{Q}_p^2}}^{G_{\mathbb{Q}_p}} \varepsilon_2^{k-1} \otimes \mu_{\sqrt{-1}}.$ Therefore,  $P_k\subset ker(\rho_{f,D}) $  and   
$$\rho_{f,D} \in Hom_{\mathbb{Q}_{p^n}}(R_{D}/P_{k}R_{D}[1/p]^{red},\bar{\mathbb{Q}}_p),$$
 where $Hom_{\mathbb{Q}_{p^n}}(R_{D}/P_{k}R_{D}[1/p]^{red},\bar{\mathbb{Q}}_p)$ is the set of ring homomorphism between $R_{D}/P_{k}R_{D}[1/p]^{red}$ and $\bar{\mathbb{Q}}_p$ which are $\mathbb{Q}_{p^n}$ linear.   Hence,  
\[m_p\left(0,k,\chi,\bar{\rho},\{\tau_v\}_{v\in S\backslash\{p\}}\right) \leq \dim_{\mathbb{Q}_{p^n}}\left( R_{D}/P_kR_{D}[1/p]^{red} \right),\]
where $ \dim_{\mathbb{Q}_{p^n}}\Big( R_{D}/P_kR_{D}[1/p]^{red} \Big)$ is the dimension of $ R_{D}/P_kR_{D}[1/p]^{red}$ as a $\mathbb{Q}_{p^n}$ vector space. 
  Assume that $R_{D}$ is finitely generated as a $\Lambda_{p}$ module with $m$ generators. Since $\Lambda_{p}/P_k=\mathbb{Z}_{p^n},$  it follows that  
$$
 m_p\left(0,k,\chi,\bar{\rho},\{\tau_v\}_{v\in S\backslash\{p\}}\right) \leq \dim_{\mathbb{Q}_{p^n}}\left( R_D/P_kR_D[1/p]^{red} \right) \leq m.
$$
Let $M$ be the maximal ideal of $\Lambda_p.$ By the topological form of the Nakayama's lemma (see \cite[Exercise 7.2]{Eisenbud}), if $R_{D} \otimes_{\Lambda_p/M} \mathbb{F}_{p^n}$ is spanned by the image of $r_1, \dots, r_{s}\in R_{D} $ over $ \mathbb{F}_{p^n},$ then $r_1, \dots, r_{s}$ generate $R_{D}$  as an $\Lambda_p$ module. Therefore,
$$
m_p\left(0,k,\chi,\bar{\rho},\{\tau_v\}_{v\in S\backslash\{p\}}\right) \leq   \dim_{\mathbb{F}_{p^n}}\left( R_D   \otimes_{\Lambda_p/M} \mathbb{F}_{p^n} \right).
$$
Recall the prime  $P_{\text{min}}$ of $\Lambda_p$ that  is associated to the dihedral representation $\text{Ind}_{G_{\mathbb{Q}_p^2}}^{G_{\mathbb{Q}_p}} \varepsilon_2^{k^{\prime}-1} \otimes \varepsilon^i\mu_{\sqrt{-1}}.$ We have 
\[
  \dim_{\mathbb{F}_{p^n}}\left( R_D   \otimes_{\Lambda_p/M} \mathbb{F}_{p^n} \right)= \dim_{\mathbb{F}_{p^n}}\Big( R_D/P_{\text{min}} R_D   \otimes_{  \Lambda_p/M} \mathbb{F}_{p^n} \Big).
\]
We give an upper bound on $\dim_{\mathbb{F}_{p^n}}\Big( R_D/P_{\text{min}} R_D   \otimes_{  \Lambda_p/M} \mathbb{F}_{p^n} \Big).$ Let $$\pi_{\min}: R_D \to R_D/P_{\min}.$$ Let $\rho_{\min}:=\pi_{\min}\circ\rho_D\otimes \varepsilon^{-i}: G_{\mathbb{Q}}  \to  GL_2(R_{D}) \to GL_2(R_D/P_{\text{min}} R_D).$

In the following lemma, we  compute  the determinant of $\rho_{\min}$. 
\begin{lem}\label{detlem}
 We have 
$$
det(\rho_{\min})=i(\chi)\varepsilon^{k^{\prime}-1}.
$$ 
\end{lem}
\begin{proof}
Recall that $\rho_D$ is unramified outside $S$ and for $v \in S$ and $v\neq p$, we have $\det(\rho_D)|_{G_{\mathbb{Q}_v}}=i(\chi_v).$ Hence,  $det(\rho_{\min})^{-1}i(\chi)$ is unramified everywhere except at prime $p.$  At $p,$ we know that the representation is isomorphic to $\text{Ind}_{G_{\mathbb{Q}_p^2}}^{G_{\mathbb{Q}_p}} \varepsilon_2^{k^{\prime}-1} \otimes \mu_{\sqrt{-1}}$ which has determinant  $\varepsilon^{k^{\prime}-1}.$  Therefore, the character $$i(\chi)\varepsilon^{k^{\prime}-1} det(\rho_{\min})^{-1}: G_{\mathbb{Q}} \to R_D/P_{\text{min}} R_D,$$
is unramified everywhere and the only unramified character on $G_{\mathbb{Q}} $ is the identity character. Hence 
$$
det(\rho_{\min})=i(\chi)\varepsilon^{k^{\prime}-1}. 
$$ 
This concludes our lemma. 
\end{proof}
 Since $k^{\prime}-1\leq \frac{p+1}{2} \leq  p-2 $ for $p\geq 5,$ by Theorem~\ref{rama}, there exists  $\bar{R}_p^{\square,\varepsilon^{k^{\prime}-1}}$ which is the unique reduced, $p$-torsion free quotient of $R_p^{\square}$ with the property that a continuous homomorphism $\alpha: R_p^{\square} \to \bar{\mathbb{Q}}_p $ factors through $\bar{R}_p^{\square,\varepsilon^{k^{\prime}-1}}$ if and only if
$\alpha \circ \rho^{\text{univ}}$ is crystalline, $\det(\alpha \circ \rho^{\text{univ}})= \varepsilon^{k^{\prime}-1}$ and   $HT(\psi \circ \rho^{\text{univ}})=(0,k^{\prime}-1)$.  Note that the dihedral representation $\text{Ind}_{G_{\mathbb{Q}_p^2}}^{G_{\mathbb{Q}_p}} \varepsilon_2^{k^{\prime}-1} \otimes \mu_{\sqrt{-1}}$ is crystalline with determinant $ \varepsilon^{k^{\prime}-1}$ and its Hodge-Tate's weights are $(0,k^{\prime}-1).$

Recall the notations while defining~\eqref{defring}. Let $\bar{R}_S^{\square,\psi}:=R_S^{\square,\psi} \hat{\otimes}_{R^{\square,loc,\psi}}\bar{R}^{\square,loc,\psi},$ where $\psi:=i(\chi)\varepsilon^{k^{\prime}-1}$ and $\bar{R}^{\square,loc,\psi}:=\bar{R}_p^{\square,\varepsilon^{k^{\prime}-1}}\hat{\otimes}_{\mathbb{Z}_{p^n},v\in {S\backslash\{p\}}}R_{v,\tau_v}^{\square,\chi_v},$ where the completed  tensor product is over $\mathbb{Z}_{p^n}$ and the inertial types $\{\tau_v\}$ are the one that we fixed for defining $m_p\left(0,k,\chi,\bar{\rho},\{\tau_v\}_{v\in S\backslash\{p\}}\right).$ Recall that $\bar{R}_S^{\psi}$ is the image of $R_{\{\tau_v,\chi_v\}_{v\in S\backslash\{p\}}}$ in $\bar{R}_S^{\square,\psi}$.

Recall that $R_{D}:=R_{\{\tau_v,\chi_v\}_{v\in S\backslash\{p\}}} \hat{\otimes}_{R_p}  \Lambda_{p}.$ The map $R_{\{\tau_v,\chi_v\}_{v\in S\backslash\{p\}}}\to R_D/P_{\text{min}} R_D$ is surjective.  Next, we show that this map factors through $\bar{R}_S^{\psi}.$  By fixing matrices  $\alpha_v\in  Ker( GL_2(R_D) \to GL_2(\mathbb{F}_{p^n}))$ for $v\in S\backslash\{p\}$, we obtain a  map from the framed deformation rings $\hat{\otimes}_{\mathbb{Z}_{p^n},v\in {S\backslash\{p\}}}R_{v,\tau_v}^{\square,\chi_v} \to R_D/P_{\text{min}} R_D$. By Lemma~\ref{detlem}, we have $det(\rho_{\min})|_{G_{\mathbb{Q}_p}}=\varepsilon^{k^{\prime}-1}$. Hence, by fixing $\alpha_p\in  Ker( GL_2(R_D) \to GL_2(\mathbb{F}_{p^n})),$ we obtain a ring homomorphism from $\bar{R}_p^{\square,\varepsilon^{k^{\prime}-1}}\to R_D/P_{\text{min}} R_D.$ 
Hence, we obtain a map  $\bar{R}_S^{\square,\psi} \to R_D/P_{\text{min}} R_D.$ Hence, we  have the following commutative diagram 
\begin{center}
\tikzcdset{
arrow style=tikz,
diagrams={>={Straight Barb[scale=0.8]}}
}
% in document body
\begin{tikzcd}
R_{\{\tau_v,\chi_v\}_{v\in S\backslash\{p\}}}\arrow[r] \arrow[rd, two heads] &\bar{R}_S^{\square,\psi} \arrow[d, two heads]\\
& R_D/P_{\text{min}} R_D.
\end{tikzcd}
\end{center}
By the above diagram the map from $\beta: \bar{R}_S^{\psi} \to R_D/P_{\text{min}} R_D$  is surjective. Hence,
\[
\dim_{\mathbb{F}_{p^n}}\Big( R_D/P_{\text{min}} R_D   \otimes_{  \Lambda_p/M} \mathbb{F}_{p^n} \Big) =  \dim_{\mathbb{F}_{p^n}}\Big( \bar{R}_S^{\psi} / ker(\beta)\otimes_{\mathbb{Z}_{p^n}} \mathbb{F}_{p^n}  \Big) .
\] 
By Theorem~\ref{finite}, $ \bar{R}_S^{\psi} $ is finite as a $\mathbb{Z}_{p^n}$ module. Hence $\dim_{\mathbb{F}_{p^n}}\Big( \bar{R}_S^{\psi}  \otimes_{\mathbb{Z}_{p^n}} \mathbb{F}_{p^n}  \Big)=O(1).$ At this point, we showed that the number of newforms $f$ (up to a scalar multiple) of level $N$ and even weight $k$ such that $\mathcal{T}_p(f)=0$  is bounded independently of $k$. We proceed and conclude Theorem~\ref{strong1.1} by showing that
\[
\dim_{\mathbb{F}_{p^n}}\Big( \bar{R}_S^{\psi} / ker(\beta)\otimes_{\mathbb{Z}_{p^n}} \mathbb{F}_{p^n}  \Big) = m_p\left(0,k^{\prime},\chi,\bar{\rho},\{\tau_v\}_{v\in S\backslash\{p\}}\right).
\]
By Corollary~\ref{maincor}, $ \bar{R}_S^{\psi}$    is finite and  Cohen-Macaulay over $\mathbb{Z}_{p^n}.$ By the discussion in \cite[Section 5]{Snowden}, it follows that $\bar{R}_S^{\psi} / ker(\beta)$ is isomorphic to a Hecke algebra $\mathbb{T}_{m}$ completed at a maximal ideal $m$ associated to the residual representation $\bar{\rho},$ where its associated local representations for $v\in S\backslash \{p\}$ has inertial type $\tau_v$ with determinant $\chi_v$ and its local representation at  $p$ is isomorphic to $\text{Ind}_{G_{\mathbb{Q}_p^2}}^{G_{\mathbb{Q}_p}} \varepsilon_2^{k^{\prime}-1} \otimes \mu_{\sqrt{-1}}.$  Therefore,  we have 
$$
\dim_{\mathbb{F}_{p^n}}\Big( \bar{R}_S^{\psi} / ker(\beta)\otimes_{\mathbb{Z}_{p^n}} \mathbb{F}_{p^n}  \Big)= \dim_{\mathbb{F}_{p^n}} (\prod_{f} \mathbb{F}_{p^{n}})
$$
where  $f\in S_{k^{\prime},\chi}(N)$ is a newform with level $N,$ inertial types $\{\tau_v\}_{v\in S\backslash \{ p\}}$ and $a_p(f)=0$. By definition the number of such  $f\in S_{k^{\prime},\chi}(N)$ is $m_p\left(0,k^{\prime},\chi,\bar{\rho},\{\tau_v\}_{v\in S\backslash\{p\}}\right)$. This concludes the proof of Theorem~\ref{strong1.1}.
\end{proof}

\subsection*{Acknowledgements}\noindent
I would like to thank Professor Frank Calegari for an insightful discussion regarding this project. Specially the key observation of imposing  the  dihedral condition at prime $p$ is due to him. I would like to thank  Professors Nigel Boston and  Richard Taylor for answering my questions. I would like to thank Professors Jordan Ellenberg and Lue Pan for their comments on the earlier version of this work.  Finally, I would like to thank the referee for  his/her comments on the earlier version of this paper.

\bibliographystyle{alpha}
\bibliography{Arxiv}
\end{document}